\documentclass[12pt]{amsart}
\usepackage{epsfig,graphicx}
\usepackage{amsmath}
\usepackage{amssymb}

\usepackage{hyperref}
\hypersetup{colorlinks=true,linkcolor=blue, linktoc=none, citecolor=blue,urlcolor=blue}

\title{Corks}
\author[AKBULUT]{Selman Akbulut} 

\thanks{}

\address{G\"{o}kova Geometry Topology Institute,  Mu\u{g}la, Turkiye}
\email{akbulut.selman@gmail.com}

\begin{document}

\begin{abstract}
Remarks relating the various notions of corks.
\end{abstract}

\maketitle

\vspace{-.2in}

I recall, my advisor Kirby kindly named the object I had constructed in \cite{a2} ``Akbulut cork". Since then, this concept went through some evolutions, resulting some confusions between its various different definitions. To clarify them, I decided to collect here mutually consistent definitions. 
  
 \vspace{.1in}
 
\noindent  Definition: A {\bf loose cork} of a smooth manifold $M$ is a pair $(W,f)$, where $W$ is a smooth contractible submanifold of $M$ with involution $f:\partial W\to \partial W$, such that cutting W out of M then regluing with $f$ changes the smooth structure of $M$.  A {\bf cork} $(W,f)$ is a loose cork, which is a Stein manifold.

 \vspace{.1in}

 We call the operation of cutting $W$ from $M$ and regluing by the map $f$ {\it a loose cork-twisting operation on $M$}, or simply cork-twisting operation on $M$. At first glance the ``Stein'' condition might look artificial, but it really is not, because by \cite{am2} every loose cork $(W,f)$ contains a cork, hence cork is a smaller fundamental object. Furthermore, given a contractible manifold with involution on its boundary $(W,f)$, it is a difficult to decide whether $W$ is a loose cork  (is there a corresponding $M$?), or does $f$ extend to a self diffeomorphism of $W$?. To answer the second question in the negative, it suffices to find a pair of knots $\gamma, \gamma' \subset \partial W $ such that $\gamma$  is not slice but $\gamma'$ is slice in $W$, with $f(\gamma)=\gamma'$. As in the example of Figure 9.1 of  \cite{a1}  (if $\gamma$ was slice you would get a violation to the ``Adjunction Inequality for Stein manifolds'', which had been proven in \cite{am1}).
 
 \vspace{.1in}
 
 To find such examples $(W,\gamma, \gamma')$, symplectic topology can be useful. For example, when $W$ is a Stein manifold, attaching $W$ a $2$-handle along a knot $\gamma\subset W$ with $TB(\gamma)-1$ framing extends the Stein structure from $W$ to $W\smile h^{2}_\gamma$ (In fact, by \cite{lm} and \cite{ao}, we can even compactify any Stein manifold into a closed symplectic manifold which could play the role of $M$ in the above definition of cork). Then by the Adjunction inequality, proved for the Stein manifolds in \cite{am1}, applied to $W\smile h^{2}_\gamma$, shows that $\gamma$ can not be slice.

 \vspace{.05in}
  
 After I wrote \cite{a2}, I realized that, even though the cork, I constructed there, looks like the manifold which Mazur had defined in his Annals paper, it is not. Mazur's manifold is not even a ``loose cork''! 
 
\vspace{.05in}

Matveyev \cite{m} generalized \cite{a2} by proving any exotic copy $M'$ of a simply connected closed smooth $4$-manifold $M$ can be obtained by performing a loose cork-twisting operation on $M$. Hence by \cite{am2} it can be obtained by a cork-twisting operation on $M$.
Curtis-Freedman-Hsiang-Stong \cite{cfhs} proved a weaker version of \cite{m} without mentioning a cork, namely  ``$M'$ can be obtained by cutting out a smooth compact contractible submanifold of $M$, then gluing another smooth compact contractible manifold with the same boundary".  
 
 \vspace{.05in}
 
  By weakening the involution condition $f$ on a loose cork $(W,f)$ one can construct an infinite order loose cork;  \cite{a3} and \cite{g} are such examples. Notice that in minute 1:14:20 of my talk
 in \cite{a4} I am trying to construct an infinite order cork (with``Stein'' condition), which we still don't know if exists.

\end{document}